\documentclass[12pt,reqno]{amsart}
\usepackage{amssymb}
\usepackage{graphicx}
\usepackage{epsf}
\usepackage{amsmath}
\usepackage[latin1]{inputenc}
\usepackage{amsfonts}
\usepackage{color}

\numberwithin{equation}{section}
\newtheorem{theo}{Theorem}[section]
\newtheorem{prop}[theo]{Proposition}
\newtheorem{conj}[theo]{Conjecture}

\newtheorem{defi}[theo]{Definition}

\newtheorem{rema}[theo]{Remark}

\newcommand{\pref}[1]{(\ref{#1})}

\newcommand{\twolines}[2]{{\scriptstyle #1 \atop \scriptstyle
\vphantom{\sum}#2}}

\def\rouge{\textcolor{red}} 
 
\def\vert{\textcolor{green}}

\def\ie{{\it i.e. }}

\def\ti{$\sim$}

\def\concat{{\cdot}}

\def\N{{\mathbb N}}

\def\Q{{\mathbb Q}}

\def\J{{\mathcal J}}

\def\DI{{\mathcal DI}}
\def\S{{\mathbb S}}

\def\B{{\mathcal B}}

\def\fh{{\mathfrak h}}

\def\QSym{{\rm QSym}}
\def\DQ{{\rm DQ}}
\def\QQ{{\rm Q}}
\def\NSym{{\rm NSym}}
\def\DNSym{{\rm DNSym}}
\def\DSSym{{\rm DSSym}}
\def\SSym{{\rm SSym}}
\def\DQSym{{\rm DQSym}}
\def\DSym{{\rm DSym}}
\def\Sym{{\rm Sym}}
\def\DSH{{\rm DTLH}}
\def\DH{{\rm DH}}
\def\TL{{\rm TL}}
\def\TLH{{\rm TLH}}
\def\Ker{{\rm Ker}}
\def\sign{{\rm sign}}
\def\DI{\J}

\def\V{{\mathcal V}}

\def\shuffle{{\,\raise
1pt\hbox{$\scriptscriptstyle\cup{\mskip-4mu}\cup$}\,}}
\def\sshuffle{{\,\hbox{$\scriptscriptstyle\cup{\mskip-4mu}\cup$}\,}}

\title[Diagonal TL invariants and Harmonics]{Diagonal Temperley-Lieb Invariants\\ and Harmonics}

\author{J.-C.~Aval}\address[Jean-Christophe Aval]{LaBRI\\ Universit\'e Bordeaux 1\\ 351 cours
de 
 la Lib\'eration\\ 33405 Talence cedex\\ FRANCE}
\email{aval@labri.fr}
\urladdr{http://www.labri.fr/\ti aval/}

 \author{F.~Bergeron} \address[F. Bergeron]{D\'epartement de Math\'ematiques\\ Universit\'e  du Qu\'ebec \`a Montr\'eal\\ Montr\'eal, Qu\'ebec, H3C 3P8, CANADA} \email{bergeron.francois@uqam.ca}

\author{N.~Bergeron}  \address[N. Bergeron]{Department of Mathematics and Statistics\\ York  University\\       To\-ron\-to, Ontario M3J 1P3\\ CANADA} \email{bergeron@mathstat.yorku.ca}  \urladdr{http://www.math.yorku.ca/bergeron}

\date{\today}

\thanks{J.C. Aval is supported by EC's IHRP Programme, Network ``ACE'', grant HPRN-CT-2001-00272.} 
 \thanks{F. Bergeron is supported in part by NSERC and  FCAR}   \thanks{N. Bergeron is supported in part by CRC, NSERC and PREA}

\begin{document} 

 \begin{abstract} 
In the context of the ring $\Q[{\bf x},{\bf y}]$, of polynomials in $2\,n$ variables ${\bf x}=x_1,\ldots,x_n$ and ${\bf y}=y_1,\ldots,y_n$, we introduce the notion of diagonally quasi-symmetric polynomials. These, also called {\em diagonal Temperley-Lieb invariants}, make possible the further introduction of the space of {\em diagonal Temperley-Lieb harmonics} and {\em diagonal Temperley-Lieb coinvariant space}. We present new results and conjectures concerning these spaces, as well as the space obtained as the quotient of the ring of diagonal Temperley-Lieb invariants by the ideal generated by constant term free diagonally symmetric invariants. We also describe how the space of diagonal Temperley-Lieb invariants affords a natural graded Hopf algebra structure, for $n$ going to $\infty$. We finally show how this last space and its graded dual Hopf algebra are related to the well known Hopf algebras of symmetric functions, quasi-symmetric functions and noncommutative symmetric functions. 

\end{abstract}   
\maketitle 
{ \parskip=0pt\footnotesize \tableofcontents}
\parskip=8pt  
\section{Introduction}\label{intro}

\subsection{General framework} \label{framework} As a warm up for of our discussion, it may be best to recall some facts about the interactions between the algebras
\begin{equation}\label{dia1} 
    \Sym_n\,\hookrightarrow\, \QSym_n\,\hookrightarrow \, \Q[ {\bf x}] \,,
\end{equation}
of {\em symmetric polynomials}, {\em quasi-symmetric polynomials}, and of polynomials in $n$ variables ${\bf x}=x_1,x_2,\ldots,x_n$. 
Some of these facts involve the natural (and well known) action of the symmetric group $\S_n$ on polynomials, defined on $p({\bf x})$ in $\Q[{\bf x}]$ as
  $$\sigma\cdot p({\bf x})=p(x_{\sigma(1)},x_{\sigma(2)},\ldots,x_{\sigma(n)}),$$
for $\sigma$ in $\S_n$. Clearly, this action extends to a faithful  action of the group algebra $\Q \S_n$ on $\Q[{\bf x}]$. Furthermore, it preserves the multiplication of polynomials,  we can thus consider the {\em algebra of symmetric polynomials}:
  $$\Sym_n\ :=\ \left\{ p({\bf x})\in\Q[{\bf x}]\ |\ \sigma\cdot p({\bf x})=p({\bf x})\right\}, $$
 which is often denoted  $\Q[{\bf x}]^{\S_n}$ by other authors. We part from this convention to stay in the spirit of rest of the notation used in this paper.
A classical result of Newton states  that $\Sym_n$ is itself a polynomial algebra. This is often made explicit as follows:
   $$\Sym_n=\Q[e_1,e_2,\ldots,e_n],$$ 
 where $e_k=e_k({\bf x})$ is the coefficient of $t^k$ in $(1+x_1t)(1+x_2t)\cdots(1+x_nt)$.

In \cite{hivert} Hivert has introduced an action of $\Q \S_n$ on polynomials which does not preserve the multiplication, but still leads to very interesting results. Before describing this action, let us first recall some notation. We denote $ {\bf x}_{\bf i}^{\bf a}$, the monomial
\begin{equation}\label{monome}
   {\bf x}_{\bf i}^{\bf a}=x_{i_1}^{a_1}x_{i_2}^{a_2}\cdots x_{i_k}^{a_k},
\end{equation}
assuming here that ${\bf a}= (a_1,a_2,\ldots,a_k)$ is a sequence of positive integers,  $a_i>0$; and that the elements $i_m$, of the subset ${\bf i}=\{i_1,i_2,\ldots ,i_k\}$ of $[n]:=\{1,\dots,n\}$,  have been listed in increasing order ($i_m<i_{m+1}$). We underline that in (\ref{monome}), the {\em support subset} ${\bf i}$ is  of cardinality equal to the length, $\ell({\bf a}):=k$, of the {\em composition exponent} ${\bf a}$. 
Recall that ${\bf a}=(a_1,a_2,\ldots,a_k)$, as above, is said to be a {\em composition} of $d$, if
        $$d=|{\bf a}|:=a_1+a_2+\ldots +a_k.$$
In our setup, $d$ is clearly the degree of  ${\bf x}_{\bf i}^{\bf a}$.
\begin{defi}[Action of Hivert]\label{hivaction}
For any permutation $\sigma$ in $\S_n$, we consider   the linear action 
    $$p({\bf x})\mapsto \sigma * p({\bf x}),$$ 
 whose effect on monomials is given as
\begin{equation}\label{hivaction_mon}
\sigma*{\bf x}_{\bf i}^{\bf a}:={\bf x}_{\sigma\cdot {\bf i}}^{\bf a},
\end{equation}
where, as usual, $\sigma\cdot {\bf i}$ is the set of $\sigma(i_m)$'s. 
\end{defi}
\noindent In the right hand side of (\ref{hivaction_mon}), one must be careful to list the elements of $\sigma\cdot {\bf i}$ in  increasing order for the purpose of expressing the corresponding monomial. For example,  applying this action for $\sigma=3\, 1\, 4\, 2$ on $x_1^2x_3x_4^2$ we get 
      $$\sigma*x_1^2x_3x_4^2= x_2^2x_3x_4^2,$$
 since $\sigma\cdot \{1,3,4\} = \{2,3,4\}$.

This action is not faithful as an algebra action. However, it can be shown (see \cite{hivert}) that the quotient of $\Q \S_n$  by the kernel of the action is precisely the Temperley-Lieb algebra $\TL_n$. Recall that this algebra can be described as follows. Consider the elements $e_i:={(q-T_i)}/{(1+q)}$ in the Hecke algebra, with $T_i$'s being the
standard generators.
The {\em Temperley-Lieb algebra}, $\TL_n(q)$, 
is defined (see \cite{jones}) as the quotient of the Hecke algebra by the relations
\begin{equation}\label{rel_TL}
   e_{i}e_{i\pm 1}e_{i} - {q\over (1+q)^2} e_{i}=0,\qquad i=1,\ldots,n
\end{equation}
The set of invariant polynomials for the Hivert action, known as {\em quasi-symmetric polynomials}, form a subalgebra of $\Q[{\bf x}]$, denoted $\QSym_n$. In formula:
  $$\QSym_n\ :=\ \left\{ p({\bf x})\in\Q[{\bf x}]\ |\ \sigma*p({\bf x})=p({\bf x})\right\}.$$
The fact that this set is closed under multiplication is not obvious from this definition. It follows from work of Stanley  \cite{stanley1} and Gessel \cite{ges}.

It is well known that the respective multiplication rules of the algebras $\Sym_n$ and $\QSym_n$ both stabilize as $n\to\infty$. It thus makes sense to consider the corresponding inverse limit algebras $\Sym$ and $\QSym$, respectively of {\em symmetric series} and {\em quasi-symmetric series}. A distinctive feature of these limit algebras is that they admit a comultiplication which is compatible with the multiplication. In this way,  both $\Sym$ and $\QSym$ actually become graded Hopf algebras. Moreover, as graded Hopf algebras,  $\Sym$  is isomorphic to its graded Hopf dual. By contrast, the graded dual of $\QSym$ is known to be another important graded Hopf algebra, namely  $\NSym$ the Hopf algebra of noncommutative symmetric functions (see \cite{MR,NC}).  The inclusion of $\Sym$ in $\QSym$ is itself a graded algebra morphism, and its dual is an epimorphism (of graded Hopf algebras)
\begin{equation}\label{diad1}   
    \NSym  \twoheadrightarrow \Sym \,.
\end{equation}
Among the many striking  other properties of $\QSym$, it is shown in \cite{ABS} that it is the terminal object in the category of {\em combinatorial Hopf algebras}. The objects of this category are graded Hopf algebras equipped with a ``character'' (multiplicative linear  functional) with value in the ground field. Arrows are morphism of graded Hopf algebras, compatible with this selected character. Any such combinatorial Hopf algebra affords two special Hopf subalgebras, called {\em even} and {\em odd}, and it is further shown in {\em ibid.} that the odd subalgebra of $\QSym$ is none other than Stembridge's \cite{Stem} {\em peak subalgebra} of quasi-symmetric functions. On another note, as shown in  \cite{NCSF4}, $\QSym$ is isomorphic to the Grothendieck ring of the projective representations of the Hecke algebra at $q=0$. It is also worth mentioning that that $\QSym$ is Free (see \cite{Hoff}) as a $\Sym$-module. The fact that this is also true for $\QSym_n$, as a $\Sym_n$-module, was first conjectured by the second author and C.~Reutenauer, and later shown to be true by A.~Garsia and N.~Wallach in \cite{GW}.

Another natural line of inquiry leads to the study of quotients associated to the algebras in (\ref{dia1}), as well as many other similar situations some of which appear in the  diagram of Figure \ref{fig_diag}. Most of the spaces of this diagram will at least be mentioned in our discussion, and all arrows have a significant role. However, it would be too long to described all these in full detail. Still, it is interesting to bear this diagram in mind as a background map for our exploration.
\begin{equation}\label{fig_diag}
{\setlength{\unitlength}{1mm}
\begin{picture}(0,30)(10,-10)
\put(-38,15){$\DNSym_n$}
\put(-30,12){\rouge{\vector(0,-1){15}}}
\put(-30,-3){\rouge{\vector(0,-1){2}}}
\put(-27,13){\vert{\vector(1,-1){9}}}
\put(-36,-10){$\DSym_n$}
\put(-27,-12){\vert{\vector(1,-1){9}}}
\put(-21,0){$\NSym_n$}
\put(-15,-3){\rouge{\vector(0,-1){15}}}
\put(-15,-18){\rouge{\vector(0,-1){2}}}
\put(-19,-25){$\Sym_n$}
\put(-20,16.45){\oval(1,1)[l]}
\put(-20,15.5){\vector(1,0){15}}
\put(-20,-8.55){\oval(1,1)[l]}
\put(-20,-9.5){\vector(1,0){15}}
\put(-5,1.45){\oval(1,1)[l]}
\put(-5,.5){\vector(1,0){15}}
\put(-5,-23.55){\oval(1,1)[l]}
\put(-5,-24.5){\vector(1,0){15}}
\put(-4,15){$\DSSym_n$}
\put(2,12){\rouge{\vector(0,-1){15}}}
\put(2,-3){\rouge{\vector(0,-1){2}}}
\put(5,13){\vert{\vector(1,-1){9}}}
\put(-4,-10){$\DQSym_n$}
\put(5,-12){\vert{\vector(1,-1){9}}}
\put(12,0){$\SSym_n$}
\put(18,-3){\rouge{\vector(0,-1){15}}}
\put(18,-18){\rouge{\vector(0,-1){2}}}
\put(12,-25){$\QSym_n$}
\put(15,16.45){\oval(1,1)[l]}
\put(15,15.5){\vector(1,0){15}}
\put(15,-8.55){\oval(1,1)[l]}
\put(15,-9.5){\vector(1,0){15}}
\put(30,1.45){\oval(1,1)[l]}
\put(30,.5){\vector(1,0){15}}
\put(30,-23.55){\oval(1,1)[l]}
\put(30,-24.5){\vector(1,0){15}}
\put(31,15){$\Q\langle{\bf x},{\bf y}\rangle$}
\put(37,12){\rouge{\vector(0,-1){15}}}
\put(37,-3){\rouge{\vector(0,-1){2}}}
\put(41,13){\vert{\vector(1,-1){9}}}
\put(31,-10){$\Q[{\bf x},{\bf y}]$}
\put(41,-12){\vert{\vector(1,-1){9}}}
\put(46,0){$\Q\langle{\bf x}\rangle$}
\put(52,-3){\rouge{\vector(0,-1){15}}}
\put(52,-18){\rouge{\vector(0,-1){2}}}
\put(46,-25){$\Q[{\bf x}]$}
\end{picture}}
\end{equation}
\vskip .6in

 The prototypal object along these lines is certainly  $\Q[{\bf x}]\big / \langle \Sym_n^+\rangle$, where 
$\langle \Sym_n^+\rangle$ denotes the ideal of $\Q[{\bf x}]$ generated by constant term free symmetric polynomials. Often encountered in association with $\Sym_n$, this quotient  plays an important role in Invariant Theory \cite{steinberg}, Galois Theory \cite{artin}, and in Algebraic Geometry \cite{borel}. In the first of these contexts, it is known as the {\em coinvariant space} of the symmetric group; and in the last, it appears as the ring of cohomology of the flag variety. It has finite dimension equal to $n!$, and it is in fact isomorphic (as a $\S_n$-module) to the left regular representation of $\S_n$. 

To describe another important feature of all this setup, we need to consider the scalar product, on $\Q[{\bf x}]$, defined by
\begin{equation}\label{scalX}
\langle p,q\rangle:=L_0\big(p(\partial {\bf x})q({\bf x})\big)
\end{equation}
where $p(\partial {\bf x})$ is the differential operator obtained by replacing  each variable $x_i$ in $p({\bf x})$, by the partial derivative $\partial x_i$. Here, $L_0$ is simply the linear form that sends a polynomial to its constant term. One easily checks that the set of monomials forms an orthogonal linear basis of $\Q[{\bf x}]$ for this scalar product.
The ideal $ \langle \Sym_n^+\rangle$ being homogeneous and $S_n$-invariant, its orthogonal complement
\begin{equation}\label{SHarm}
     \mathcal{H}_n :=  \langle \Sym_n^+\rangle^\perp.
\end{equation}
is a graded $\S_n$-module (by degree). It is easy to check that $q({\bf x})$ is in $\mathcal{H}_n$ if and only if it satisfies all partial differential equation of the form 
  $$p(\partial {\bf x}) q({\bf x})=0,$$
with $p({\bf x})$ any constant term free symmetric polynomial.
In particular, observing that $x_1^2+\cdots+x_n^2$ is evidently symmetric, it follows that not only are all elements of $\mathcal{H}_n$ harmonic polynomials in the classical sense, but they further satisfy several other conditions:
   $$\partial x_1^k \,p({\bf x}) + \dots + \partial x_n^k\, p({\bf x})=0, \qquad 1\leq k\leq n. $$
 For this reason, the space $\mathcal{H}_n$ is know as the space of $\S_n$-harmonic polynomials. In fact, $\mathcal{H}_n$ can be explicitly described as being  simply the linear span of all partial derivative of the Vandermonde determinant:
     $$\Delta_n({\bf x}):=\prod_{i<j} (x_i-x_j).$$
Now, for the same reason as for $\mathcal{H}_n$, the space $\Q[{\bf x}]\big / \langle \Sym_n^+\rangle$ is graded, and the two spaces are actually isomorphic as graded $S_n$-modules. In truth, all of these nice properties stem from the fact that $\S_n$ acts as a reflection group on the linear part of $\Q[{\bf x}]$. This is because adapted statements of this kind actually characterize finite reflection groups (see \cite{steinberg}).

In a recent paper \cite{ABBQ} we have considered the quotient 
\begin{equation}\label{abbq}
   \QQ_n:=\Q[{\bf x}]\big / \langle \QSym_n^+\rangle.
 \end{equation}
  Surprisingly, the dimension of this space is exactly the $n^{\rm th}$ Catalan number, which correspond to the dimension of $\TL_n$. This suggests that, in analogy with the construction above, we consider the space
\begin{equation}\label{SHarm}
   \TLH_n :=  \langle \QSym_n^+\rangle^\perp \subset \mathcal{H}_n,
\end{equation}
which is clearly isomorphic (as a graded module) to $\QQ_n$.
In view of the discussion above, it seems natural to say that this is the space of {\em  $\TL_n$-harmonics}. Besides  the results of \cite{ABBQ}, not much is know about this space and it should certainly be investigated further. In particular, it would be nice to have as simple a description of it as that of $\mathcal{H}_n$.
On a similar note, and assuming the freeness of $\QSym_n$ over $\Sym_n$ (which has been shown to be true in \cite{GW}, as we have mentioned above),  it is easy to check that $\QSym_n / \langle \Sym_n^+\rangle$ is of dimension $n!$. This was originally observe by the second author and C.~Reutenauer when they conjectured this freeness property. On top of deducing an explicit expression for the Hilbert series of this same space, they further conjectured an explicit description of a permutation indexed basis (see (\ref{qsym_sym})). This basis conjecture is still open.

\subsection{Our context} The purpose of our endeavor is to extend all of the considerations above to the ``diagonal context''. Studies of the diagonal action of the symmetric group on polynomials has led to very surprising results and applications, see
\cite{lattice, nabla, GH, orbit,haiman} for a start. Our current work proposes an analogous  diagonal theory for $\TL_n$-invariants. This is a long term project, and we will not have the space to include all possible avenues. We will limit ourselves to the foundation and the ground breaking results.
In particular, we will only consider the algebras of the following diagram, which is but the bottom part of the diagram appearing in (\ref{fig_diag}).
\begin{equation}\label{dia2} {\setlength{\unitlength}{1mm} \begin{picture}(40,15)(0,0) 

\put(-21,10){$\DSym_n$}
\put(-15,7){\rouge{\vector(0,-1){15}}}
\put(-15,-8){\rouge{\vector(0,-1){2}}}
\put(-19,-15){$\Sym_n$}

\put(-5,11.45){\oval(1,1)[l]}
\put(-5,10.5){\vector(1,0){15}}
\put(-5,-13.55){\oval(1,1)[l]}
\put(-5,-14.5){\vector(1,0){15}}

\put(12,10){$\DQSym_n$}
\put(18,7){\rouge{\vector(0,-1){15}}}
\put(18,-8){\rouge{\vector(0,-1){2}}}
\put(12,-15){$\QSym_n$}

\put(30,11.45){\oval(1,1)[l]}
\put(30,10.5){\vector(1,0){15}}
\put(30,-13.55){\oval(1,1)[l]}
\put(30,-14.5){\vector(1,0){15}}

\put(48,10){$\Q[ {\bf x},{\bf y}]$}
\put(54,7){\rouge{\vector(0,-1){15}}}
\put(54,-8){\rouge{\vector(0,-1){2}}}
\put(48,-15){$\Q[ {\bf x}]$}

\end{picture}}
\end{equation}

\vskip .65in
\noindent
All of these algebras sit inside the algebra $\Q[{\bf x},{\bf y}]$ of polynomials in two sets of $n$ variables, with ${\bf y}=y_1,y_2,\ldots,y_n$.
On this space we consider the natural diagonal action of the symmetric group $\S_n$ on $\Q[{\bf x},{\bf y}]$ which permutes variables simultaneously
 $$\sigma\cdot x_i=x_{\sigma(i)},\qquad {\rm and}\qquad 
     \sigma\cdot y_i=y_{\sigma(i)}.$$
This results in a faithful  action of the group algebra $\Q \S_n$ on $\Q[{\bf x},{\bf y}]$, which clearly preserves the multiplication of polynomials. Hence, the associated polynomial invariants form an algebra, called the {\em algebra of diagonally symmetric polynomials}:
  $$\DSym_n\ :=\ \left\{ p({\bf x},{\bf y})\in\Q[{\bf x},{\bf y}]\ |\ \sigma\cdot p({\bf x},{\bf y})=p({\bf x},{\bf y})\right\} .  $$
We naturally extend the scalar product  (\ref{scalX}) to $\Q[{\bf x},{\bf y}]$, setting
\begin{equation}\label{scalXY}
\langle p,q\rangle:=L_0\big(p(\partial {\bf x},\partial {\bf y})q({\bf x}, {\bf y})\big).
\end{equation}
The corresponding space of {\em $\S_n$-diagonal harmonics}, is then defined as
\begin{equation}\label{SDHarm}
   \DH_n :=  \langle \DSym_n^+\rangle^\perp\  .
\end{equation}
It is isomorphic, as a bigraded $\S_n$-module, to the {\em diagonal coinvariant space}:
    $$\Q[{\bf x},{\bf y}] \big / \langle \DSym_n^+\rangle.$$
In one of these two incarnations, this space has been studied extensively  in the last 15 years by several authors (see for example  \cite{lattice, nabla, GH, orbit,haiman}). To finish the introduction of the spaces appearing in (\ref{dia2}), we diagonally extend to $\Q[{\bf x},{\bf y}]$ the action of $\Q \S_n$ characterized by (\ref{hivaction_mon}). This gives a faithful action of $\TL_n$ on $\Q[{\bf x},{\bf y}]$ which enables the introduction of the space $\DQSym_n $ of diagonally quasi-symmetric polynomials. All of this will be described in more details in the next section.

This paper is divided into three main sections.  First, we study the diagonal Temperley-Lieb polynomial invariants. Second, we develop the Hopf algebra theory behind this object. This includes  describing the graded dual Hopf algebra as well as showing its freeness and cofreeness. 
Let us point out that that S.~Poirier \cite{Poirier} was the  first to introduce the space $\DQSym$, although in
a completely different context, but he did not study its properties, nor did he describe its
dual. More Hopf properties of these algebras and their generalization are
studied in \cite{Hohlweg, Kape, NovThi}.
Our third section is a collection of basic results and conjectures regarding these new spaces and their interplay. In particular, we discuss the structure of the space diagonal Temperley-Lieb harmonics, as well as the quotient of the ring of diagonal Temperley-Lieb invariants by the ideal generated by constant term free diagonally symmetric invariants. 

\section{$\DQSym$ and a faithful action of $\TL_n$ on $\Q[{\bf x},{\bf y}]$}
In order to setup our discussion, we need some more definitions.
Let both ${\bf a}$ and ${\bf b}$ be length $k$ vectors of non-negative integers (in $\N^k$). The pair $({\bf a},{\bf b})$ is said to be a {\em bipartite vector}, of {\em length} $\ell({\bf a},{\bf b}):=k$, for which we  use the $2\times k$ {\em matrix notation}
\begin{equation}\label{bipartitea}
     ({\bf a},{\bf b})= \begin{pmatrix} a_1 & a_2 & \ldots & a_n\\ 
                       b_1 & b_2 & \ldots & b_n \end{pmatrix}
\end{equation}
We say that  $({\bf a},{\bf b})$  is a {\em bicomposition}, if none of the $(a_i,b_i)$'s is the zero vector $(0,0)$. The {\em empty bicomposition}, of length 0, is denoted ${0 \choose 0}$.
Clearly, (\ref{bipartitea}) is equivalent to saying that a bicomposition ${\mathfrak a} =({\bf a},{\bf b})$  is a sequence      
\begin{equation}\label{bipartiteb}
   {\mathfrak a} =(\alpha_1,\alpha_2,\ldots,\alpha_k)
\end{equation}
of {\em bivectors} $\alpha_i= \big({a_i\atop b_i}\big)  \in \N\times\N$, none of which equals   $\big({0 \atop 0}\big)$. The two points of view have respective merits, hence we will use both.
The {\em bidegree}, $d({\mathfrak a} )$, of a bicomposition ${\mathfrak a} =({\bf a},{\bf b})$ is simply
\begin{equation}\label{bidegree}
  d({\bf a},{\bf b}):=(|{\bf a}|,|{\bf b}|).
 \end{equation} 
 We finally denote by $c({\bf a},{\bf b})$ the (ordinary) composition simply obtained as the vector sum ${\bf a}+{\bf b}$.

As for monomials in one set of variables (see (\ref{monome})), we use  the following notation for monomials in the two set of variables ${\bf x}$ and ${\bf y}$, setting:
\begin{eqnarray}
    ({\bf x}{\bf y})_{\bf i}^{\mathfrak a}
   &=& (x_{i_1}^{a_1}x_{i_2}^{a_2}\cdots x_{i_k}^{a_k})( y_{i_1}^{b_1}y_{i_2}^{b_2}\cdots y_{i_k}^{b_k})\label{bimonomea}\\
   &=&  (x_{i_1}^{a_1}y_{i_1}^{b_1})( x_{i_2}^{a_2}y_{i_2}^{b_2})\cdots    
             (x_{i_k}^{a_k}y_{i_k}^{b_k}).\label{bimonomeb}
\end{eqnarray}
As before, we insist here that the elements $i_m$, of the set ${\bf i}$, come in increasing order
   $$1\le i_1<i_2<\cdots <i_k\le n.$$
 In an evident manner, the two approaches to bicompositions, corresponding to (\ref{bipartitea}) and  (\ref{bipartiteb}), give the respective expressions (\ref{bimonomea}) and  (\ref{bimonomeb}).  
Just as before, a monomial is completely characterized by its support set ${\bf i}$ and its {\em bicomposition exponent} ${\mathfrak a}$, both having the same length. Moreover, the projection map 
   \begin{equation}\label{projection}
      \pi:=\Q[{\bf x},{\bf y}]\twoheadrightarrow \Q[{\bf x}],
   \end{equation}
  which identifies $y_i=x_i$, can evidently be expressed, in term of monomials, as 
  $$ {({\bf x}{\bf y})}_{\bf i}^{\mathfrak a} \mapsto  {\bf x}_{\bf i}^{c({\mathfrak a})}.$$
  This is one of the maps involved in both diagrams (\ref{dia1}) and (\ref{dia2}).

\subsection{Action of $\TL_n$.} We can now define our diagonal extension of Hivert's action to $\Q[{\bf x},{\bf y}]$ by its effect on monomials:
\begin{equation}\label{Dhivaction}
   \sigma* {({\bf x}{\bf y})}_{\bf i}^{\mathfrak a}= {({\bf x}{\bf y})}_{\sigma\cdot {\bf i}}^{\mathfrak a}.
\end{equation}
Clearly, this extension is compatible with the projection map from $\Q[{\bf x},{\bf y}]$ to $\Q[{\bf x}]$.
For example, with $\sigma=3\, 1\, 4\, 2$, we get
   $$\sigma*x_1^2y_1\,y_3\,x_4^2= x_2^2y_2\,y_3\,x_4^2,$$
and the compatibility with $\pi$ is apparent in
\begin{eqnarray*}
    \sigma*\pi(x_1^2y_1\,y_3\,x_4^2)&=&  \sigma*x_1^3\,x_3\,x_4^2\\
       &=& x_2^3\,x_3\,x_4^2\\
       &=& \pi(x_2^2y_2\,y_3\,x_4^2)\\
       &=&\pi(\sigma*x_1^2y_1\,y_3\,x_4^2)
\end{eqnarray*}
Just as in \cite{hivert}, we have the following proposition. Observe here that the relation (\ref{rel_action}), appearing below, is the specialization at $q=1$ of the defining relation (\ref{rel_TL}) for the Temperley-Lieb algebra inside the Hecke algebra (for i=1).

\begin{prop}
For $n\ge 3$, the kernel $\Ker_n\subseteq \Q\S_n$ of the action defined by~(\ref{Dhivaction}), is the two-sided ideal generated by the element
\begin{eqnarray}\label{rel_action}
E_{\{1,2,3\}} & = & 1\,2\,3\,4\cdots n\ -\ 1\,3\,2\,4\cdots n\ -\ 2\,1\,3\,4\cdots n\nonumber\\
&&+\ 2\,3\,1\,4\cdots n\ +\ 3\,1\,2\,4\cdots n\ -\ 3\,2\,1\,4\cdots n.\label{rel_action}
\end{eqnarray}
Thus, $\TL_n\simeq \Q\S_n\big/ \Ker_n$ acts faithfully on $ \Q[{\bf x},{\bf y}]$.
\end{prop}

\begin{proof} This is just an adaptation of the proof in \cite{hivert}.  Given any monomial $({\bf x}{\bf y})_{\bf i}^{\mathfrak a}$, we first check (case by case) that $E_{\{1,2,3\}}*({\bf x}{\bf y})_{\bf i}^{\mathfrak a}=0$, or in other words that $E_{\{1,2,3\}}$ is in $\Ker_n$.
By definition. $E_{\{1,2,3\}}$ is the signed sum of the permutation that fix all but the the numbers in $\{1,2,3\}$. The various cases we have to check correspond to the number of elements in ${\bf i}\cap\{1,2,3\}$.
\begin{enumerate}
 \item[1)] If ${\bf i}\cap\{1,2,3\}$ equals either $\{1,2,3\}$ or $\emptyset$, then each of six permutations $\sigma$ in the support of $E_{\{1,2,3\}}$ is such that $\sigma*({\bf x}{\bf y})_{\bf i}^{\mathfrak a}=({\bf x}{\bf y})_{\bf i}^{\mathfrak a}$, and clearly 
\begin{equation}\label{la_rel}
    E_{\{1,2,3\}}*({\bf x}{\bf y})_{\bf i}^{\mathfrak a}=0.
 \end{equation}
 \item[2)] If  the number of elements in ${\bf i}\cap\{1,2,3\}$ is $1$ or $2$, then for any ${\bf j}$ such that $\sigma*({\bf x}{\bf y})_{\bf i}^{\mathfrak a}=({\bf x}{\bf y})_{\bf j}^{\mathfrak a}$ there are exactly two permutations $\sigma,\sigma'$ in the support of $E_{\{1,2,3\}}$ such that 
   $$\sigma*({\bf x}{\bf y})_{\bf i}^{\mathfrak a}=({\bf x}{\bf y})_{\bf j}^{\mathfrak a}=\sigma'*({\bf x}{\bf y})_{\bf i}^{\mathfrak a}.$$
 Moreover, these two permutations are such that  $\sign(\sigma)=-\sign(\sigma')$, implying that (\ref{la_rel}) holds for these cases also.
 \end{enumerate}

\noindent To finish the argument, we first observe that, by conjugation, we have that $E_{\{i,j,k\}}$ is in the two-sided ideal generated by $E_{\{1,2,3\}}$, for any $1\le i<j<k\le n$. Here $E_{\{i,j,k\}}$ is, as before, the signed sum of the permutations of that fix all but the numbers in  $\{i,j,k\}$. On the other hand, the quotient of $\Q\S_n$ by the two-sided ideal generated by the $ E_{\{i,j,k\}}$'s (or simply $E_{\{1,2,3\}}$ in view of the above argument) is well known to be $\TL_n$ (see for example \cite{stanley2}).  We have shown that in the first part of the proof that this ideal $\langle E_{\{1,2,3\}}\rangle$ is contained in $\Ker_n$. Since $\Q\S_n$ is semi-simple, the quotient $\TL_n$ of $\Q\S_n$ by a two sided ideal is also semi-simple. To show the other inclusion, it is thus enough to show that each irreducible representations of $\TL_n$ occurs in $\Q[{\bf x},{\bf y}]$. To this end, for any bicomposition ${\mathfrak a}$, consider  the invariant subspace $K_{\mathfrak a}$, of $\Q[{\bf x},{\bf y}]$,
spanned by all monomials $({\bf x}{\bf y})_{\bf i}^{\mathfrak a}$, with  ${\bf i}$ of cardinality equal to length of ${\mathfrak a}$.
It is readily seen that the corresponding representation of $\S_n$, on this invariant subspace $K_{\mathfrak a}$, is equivalent to the induced representation from $\S_k\times\S_{n-k}$ to $\S_n$. It is also well known that this representation contains all and only the irreducible $V_\lambda$ indexed\footnote{With the usual canonical indexing of irreducible representations of $S_n$ by partitions.} by partitions $\lambda$ with no more than two parts. This exhausts all irreducible representation of $\TL_n$ (see \cite{stanley2}). We conclude that $\TL_n$ acts faithfully on $\Q[{\bf x},{\bf y}]/\langle E_{\{1,2,3\}}\rangle$.
 \end{proof}
 
 \begin{proof}[\bf Remark 1.] It is interesting to describe the full bigraded character of $\Q[{\bf x},{\bf y}]$ as a  $\TL_n$-representation. Recall that the {\em Frobenius characteristic} is defined on a representation $\mathcal{W}$, of $\S_n$, as
     $$  \mathcal{F}(\mathcal{W}):=\frac{1}{n!}\sum_{\sigma\in \S_n} \chi^\mathcal{W}(\sigma)\, p_{\lambda(\sigma)},$$
where $\chi^\mathcal{W}$ denotes the character of $\mathcal{W}$, with $\lambda(\sigma)$ being the partition describing the cycle lengths in $\sigma$, and $p_\lambda$ denoting the usual power sum symmetric functions. This map sends irreducible representation indexed by $\lambda$ to the Schur symmetric function $s_\lambda$ (in the natural indexing).
 Given a bicomposition ${\mathfrak a}$,  of length $k=\ell({\mathfrak a})$, the Frobenius characteristic of the space $K_{\mathfrak a}$ (as defined in the proof of the proposition) is
     $${\mathcal F}(K_{\mathfrak a})=h_kh_{n-k},$$ 
  where $h_k$ is the usual $k^{\rm th}$ complete homogeneous symmetric function. The {\em bigraded Frobenius characteristic} of $\Q[{\bf x},{\bf y}]$, as a $\TL_n$-module, is thus given by
  $${\mathcal F}_{(q,t)}(\Q[{\bf x},{\bf y}]) = \sum_{k\ge 0} \left( {q+t-qt \over (1-q)(1-t)} \right)^k
      h_kh_{n-k}.$$
This meaning of ``bigraded'', in the name above, is that the coefficient of $q^at^b$ in this last formula is the Frobenius characteristic of the homogeneous component of $\Q[{\bf x},{\bf y}]$ of bidegree $(a,b)$.
\end{proof}

\subsection{The algebra $\DQSym_n$.} We can now define the ``algebra'' of {\em diagonally quasi-symmetric} polynomials
  $$\DQSym_n\ :=\ \left\{ p({\bf x},{\bf y})\in\Q[{\bf x},{\bf y}]\ |\ \sigma*p({\bf x},{\bf y})=p({\bf x},{\bf y})\right\}.$$
Clearly, $\DQSym_n$ is a vector space with basis given by the orbit sums of monomials. More precisely, for any bicomposition ${\mathfrak a}$, let us define the {\em bimonomial} quasi-symmetric polynomial:
  $$M_{\mathfrak a}({\bf x},{\bf y}):=\sum_{{\bf i}\subseteq [n],\  |{\bf i}|=\ell({\mathfrak a})} ({\bf x}{\bf y})_{\bf i}^{\mathfrak a}.$$
Evidently this definition implies that $M_{\mathfrak a}=0$ whenever $\ell({\mathfrak a})>n$.   The set of $M_{\mathfrak a}$, for $\ell({\mathfrak a})\leq n$,  forms a linear basis of $\DQSym_n$.
This basis makes evident that the algebra $\QSym_n$  (see (\ref{abbq})) appears, up to isomorphism, as a subspaces of $\DQSym_n$ in two obvious ways. First as the linear span of the set
    $$\{\ M_{({\bf a},0)}\ |\ \ell({\bf a})\leq  n\ \},\qquad{\rm second\ as \ that\ of\ }\qquad \{\ M_{(0,{\bf b})}\ |\ \ell({\bf b})\leq  n\ \},$$
with ${\bf a}$ and ${\bf b}$ compositions. To distinguish these two occurrences, we respectively denote them $\QSym_n^{\bf x}$ and  $\QSym_n^{\bf y}$. As we have already mentioned, these are respectively subalgebras of $\Q[{\bf x}]$ and $\Q[{\bf y}]$.
To show that, in turn, $\DQSym_n$ is itself a subalgebra of $\Q[{\bf x},{\bf y}]$, we easily adapt the method of \cite{Hoff}. To this end, we introduce the sum, concatenation and semi-shuffle of bicompositions.  For ${\mathfrak a}=(\alpha_1,\ldots,\alpha_k)$ and ${\mathfrak b}=(\beta_1,\ldots,\beta_\ell)$, the {\em sum} of ${\mathfrak a}$ and ${\mathfrak b}$ is the length $k+\ell-1$ bicomposition
\begin{equation}\label{concat}
   {\mathfrak a}+ {\mathfrak b}:=(\alpha_1,\ldots,\alpha_k+\beta_1,\ldots,\beta_\ell),
\end{equation}
with $\alpha_k+\beta_1$ being the vector sum;
and the {\em concatenation} of ${\mathfrak a}$ and ${\mathfrak b}$ is the length $k+\ell$  bicomposition
\begin{equation}\label{concat}
   {\mathfrak a}\concat  {\mathfrak b}:=(\alpha_1,\ldots,\alpha_k,\beta_1,\ldots,\beta_\ell)
\end{equation}
The {\em quasi-shuffle} of bicompositions can then be defined as the set described recursively as follows. 
\begin{equation}
  {\mathfrak a} \tilde\shuffle {\mathfrak b}:=\begin{cases}
     \{\mathfrak{a}\}& \text{if\ } \mathfrak{b}=0, \\
     \{\mathfrak{b}\}& \text{if\ } \mathfrak{a}=0,\\
    \displaystyle \bigcup_{\mathfrak{a}_1\concat  \mathfrak{a}_2=\mathfrak{a}} \ 
    \bigcup_{\mathfrak{c}\in \mathfrak{a}_2\tilde\sshuffle\mathfrak{b}_2} 
        \{ \mathfrak{a}_1\concat  \mathfrak{b}_1\concat  \mathfrak{c},
           (\mathfrak{a}_1+\mathfrak{b}_1)\concat  \mathfrak{c} \}& \text{otherwise}.
\end{cases}
\end{equation}
where $\mathfrak{b}=\mathfrak{b}_1\concat  \mathfrak{b}_2$ is the unique factorization of $\mathfrak{b}$ such that $\mathfrak{b}_1$ has only one part.
For example 
 $$\begin{array}{rcl}
 {\scriptstyle \big({2\ 0 \atop 1\ 3}\big)} \tilde\shuffle {\scriptstyle \big({0\ 2 \atop 1\ 0}\big)} &=&
    \Big\{
       {\scriptstyle \big({2\ 0\ 0\ 2 \atop 1\ 3\ 1\ 0}\big)} , {\scriptstyle \big({2\ 0\ 2 \atop 1\ 4\ 0}\big)},
       {\scriptstyle \big({2\ 0\ 0\ 2 \atop 1\ 1\ 3\ 0}\big)} , {\scriptstyle \big({2\ 0\ 2 \atop 2\ 3\ 0}\big)} 
                        ,{\scriptstyle \big({2\ 0\ 2 \atop 1\ 1\ 3}\big)} ,{\scriptstyle \big({2\ 2 \atop 2\ 3}\big)} ,
       {\scriptstyle \big({2\ 0\ 2\ 0 \atop 1\ 1\ 0\ 3}\big)} ,\cr
     &&\quad   {\scriptstyle \big({2\ 2\ 0 \atop 2\ 0\ 3}\big)} ,{\scriptstyle \big({0\ 2\ 0\ 2 \atop 1\ 1\ 3\ 0}\big)} ,
                     {\scriptstyle \big({0\ 2\ 2 \atop 1\ 1\ 3}\big)} ,
       {\scriptstyle \big({0\ 2\ 2\ 0 \atop 1\ 1\ 0\ 3}\big)} ,   {\scriptstyle \big({0\ 4\ 0 \atop 1\ 1\ 3}\big)} ,
       {\scriptstyle \big({0\ 2\ 2\ 0 \atop 1\ 0\ 1\ 3}\big)} 
    \Big\}\end{array}
$$
Following the argument of \cite{Hoff}, we get
\begin{equation}\label{DMmult}
  M_{\mathfrak a}({\bf x},{\bf y}) M_{\mathfrak b}({\bf x},{\bf y}) =\sum_{{\mathfrak c}\,\in\,{\mathfrak a}\tilde\sshuffle{\mathfrak b}} M_{\mathfrak c}({\bf x},{\bf y})\,. 
  \end{equation}
  Observe that some of the terms, $M_{\mathfrak c}$, in the resulting sum  may vanish when $\ell({\mathfrak c})>n$. One notes that the resulting multiplication rule  stabilizes as $n\to\infty$, as $n$ becomes larger then the maximum possible length for ${\mathfrak c}$. Hence we can consider the corresponding inverse limit algebra, $\DQSym$, of {\em diagonally quasi-symmetric series}. As a vector space this algebra is clearly generated by the bimonomial basis elements
  $$M_{\mathfrak a}=\sum_{{\bf i}\subseteq \N,\  |{\bf i}|=\ell({\mathfrak a})} ({\bf x}{\bf y})_{\bf i}^{\mathfrak a},$$
  where  ${\mathfrak a}$ runs over all bicompositions. The distinctive feature here is that ${\bf i}$ can be any $k$-subset of $\N$, hence $M_{\mathfrak a}$ is a series rather than a polynomial. 
 The multiplication of the $M_{\mathfrak a}$'s in $\DQSym$ is just as in (\ref{DMmult}), with no vanishing terms. To introduce a graded bialgebra structure on $\DQSym$, we define the comultiplication
    $$\Delta(M_{\mathfrak{a}}) =\sum_{\mathfrak{a}_1\concat  \mathfrak{a}_2=\mathfrak{a}} M_{\mathfrak{a}_1} \otimes M_{\mathfrak{a}_2},$$
and counit 
   $$\epsilon(M_{\mathfrak a})=\begin{cases}
      1& \text{if\ } {\mathfrak a}=0, \\
      0 & \text{otherwise}.
\end{cases}$$
The fact that we actually get a bigraded Hopf algebra, with respect to bidegree, follows from general principles (see  \cite{MM65}). Here the bidegree of $M_{\mathfrak a}$ is naturally set to be the bidegree $d({\mathfrak a})$ of ${\mathfrak a}$ (see (\ref{bidegree})).
The restriction to $\DQSym_n$ of the map $\pi$, introduced in (\ref{projection}), clearly factors through  $\QSym_n$, since for every bicomposition ${\mathfrak a}=({\bf a},{\bf b})$ we evidently have
    $$\pi(M_{({\bf a},{\bf b})})=M_{{\bf a}+{\bf b}}.$$
This map preserves multiplication, and when $n\to\infty$ it also preserves comultiplication and counit.    Observe that $\DQSym$  can also be considered as a simply graded Hopf algebra by setting 
    $$\DQSym^{(n)}:=\bigoplus_{n=d_1+d_2} \DQSym^{(d_1,d_2)},$$
 with $\DQSym^{(d_1,d_2)}$ being the bihomogeneous component of bidegree $(d_1,d_2)$ of $\DQSym$. 
We get the following theorem.

\begin{theo} $\DQSym$ is a bigraded Hopf algebra and the map $\pi\colon \DQSym \to \QSym$ is a morphism of graded Hopf algebras.\qed \end{theo}

We end this section with the description of another important basis of $\DQSym$. Consider the partial order on bicompositions obtained as the transitive closure of the set of covering relations
     $$\{ {\mathfrak a}_1+{\mathfrak a}_2 < {\mathfrak a} \ |\ {\mathfrak a}_1\concat {\mathfrak a}_2 = {\mathfrak a}\}\,.$$
We introduce the following diagonal analogue of the usual ``fundamental basis'' for quasi-symmetric functions:
\begin{equation}\label{fondamental}
     F_{\mathfrak b}:=\sum_{{\mathfrak b}\ge{\mathfrak a}} M_{\mathfrak a}.
 \end{equation}
For instance, we have
 $$\begin{array}{rcl}
 F_{ \big({2\ 0 \atop 1\ 1}\big)} &=& M_{ \big({2\ 0 \atop 1\ 1}\big)} + M_{ \big({1\ 1\ 0 \atop 1\ 0\ 1}\big)} 
      + M_{ \big({1\ 1\ 0 \atop 0\ 1\ 1}\big)} + M_{ \big({2\ 0\ 0 \atop 0\ 1\ 1}\big)} 
      + M_{ \big({0\ 2\ 0 \atop 1\ 0\ 1}\big)} + \cr  \cr&&M_{ \big({1\ 1\ 0\ 0 \atop 0\ 0\ 1\ 1}\big)}
      + M_{ \big({1\ 0\ 1\ 0 \atop 0\ 1\ 0\ 1}\big)}+ M_{ \big({0\ 1\ 1\ 0 \atop 1\ 0\ 0\ 1}\big)}
      \,.
      \end{array}
$$

\section{The graded Hopf dual \DNSym, and other Hopf properties}

Consider the non-commutative free algebra 
 $$\DNSym:=\Q\langle {\mathfrak h}_\alpha \ |\ \alpha\in \N^2\setminus {0} \rangle,$$
of non commutative polynomials in the variables ${\mathfrak h}_\alpha$, as a bigraded algebra for which we set the bidegree of the variable $\fh_\alpha$,  to be  $\alpha=(a,b)$.  Here, as usual, we denote $0$ the zero vector in $\N^2$. We further introduce a bigraded coalgebra structure on $\DNSym$, setting as comultiplication 
  $$\Delta(\fh_\alpha):=\sum_{\alpha=\beta+\gamma} \fh_\beta\otimes \fh_\gamma,$$
where $\beta+\gamma$ is understood to be the usual vector sum. By convention we let $\fh_{0}:=1$ be the unity of $\DNSym$. The counity is simply the linear form such that
   $$\epsilon(\fh_\alpha):=\begin{cases}
     1 & \text{if\  }\alpha =0, \\
     0 & \text{otherwise}.
\end{cases}$$ 
It is then straightforward to check that $\DNSym$ is in fact a bigraded Hopf algebra. This bigrading corresponding to the decomposition 
   $$\DNSym=\bigoplus_{(d_1,d_2)} \DNSym_n^{(d_1,d_2)}$$
in bihomogeneous components $\DNSym_n^{(d_1,d_2)}$. Each of these bihomogeneous components afford a linear basis 
$$ \Big\{ \fh_{\mathfrak a}  \ | \ {\mathfrak a}=(\alpha_1,\cdots,\alpha_k),\quad \alpha_i=(a_i,b_i), \ {\rm and}\  d({\mathfrak a})=(d_1,d_2)\Big\}$$

\begin{theo} $\DNSym$ is isomorphic to the bigraded Hopf dual algebra $\DQSym^*$. \end{theo}

\proof As a vector space, we clearly must have 
  $$\DQSym^* =\bigoplus_{(d_1,d_2)} {\DQSym_n^{(d_1,d_2)}}^*.$$ 
 We consider the pairing
    $$\langle \fh_{\mathfrak a},M_{\mathfrak b}\rangle:=\begin{cases}
      1 & \text{if\ }{\mathfrak a}={\mathfrak b}, \\
      0 & \text{otherwise},
\end{cases}$$
which gives a graded vector space isomorphism between $\DNSym$ and $\DQSym^*$.  To complete our proof, we need to show that the multiplicative and comultiplicative structure is preserved. Let $c_{{\mathfrak a},{\mathfrak b}}^{\mathfrak c}$ be the multiplicative structure constants of $\DQSym$:
  $$M_{\mathfrak a} M_{\mathfrak b} =\sum_{{\mathfrak c}}  c_{{\mathfrak a},{\mathfrak b}}^{\mathfrak c} M_{\mathfrak c}.$$
From (\ref{DMmult}), we have 
  $$c_{{\mathfrak a},{\mathfrak b}}^{\mathfrak c} = \left\{  \begin{array}{ll}    1 &\mbox{ if } {\mathfrak c} \in {\mathfrak a}\tilde\sshuffle {\mathfrak b} , \medskip     \\      0& \mbox{ otherwise.}\end{array}
   \right.
     $$
 On the other hand, in $\DNSym$, we have 
 $$\Delta(\fh_{\mathfrak c})=\sum_{{\mathfrak c} \in {\mathfrak a}\tilde\sshuffle {\mathfrak b} } \fh_{\mathfrak a}\otimes\fh_{\mathfrak b}.$$
 This implies that 
 $\langle \fh_{\mathfrak c}, M_{\mathfrak a} M_{\mathfrak b}\rangle = \langle\Delta( \fh_{\mathfrak c}), M_{\mathfrak a}\otimes M_{\mathfrak b}\rangle$. In the same way one can show that 
 $\langle \fh_{\mathfrak a} \fh_{\mathfrak b}, M_{\mathfrak c} \rangle = \langle  \fh_{\mathfrak a}\otimes \fh_{\mathfrak b}, \Delta(M_{\mathfrak c}) \rangle$, thus ending the proof.
\qed

\medskip
A consequence of this is that $\DQSym$ is cofree, since $\DNSym$ is free.  We can also show that $\DQSym$ is free. This is a consequence of a more general fact about quasi-shuffle Hopf algebra (see \cite{Hoff}).
To do this, we consider bicompositions as words on the alphabet of biletters $\left( a \atop b\right)$ ordered lexicographically. As is usual in the context of words, we say that a bicomposition ${\mathfrak a}$ is a {\sl Lyndon} bicomposition, if it is strictly smaller than all of the non-trivial circular rearrangement of its letters. The following theorem is a special case of a theorem of \cite{Hoff}:

\begin{prop} [Hoffman] \label{hoffman}
$\DQSym=\Q[M_{\mathfrak a} \ | \ {\mathfrak a} \hbox{ is Lyndon}]$, and in particular it is free.
\end{prop}

\begin{proof}[\bf Remark 2.] S. Poirier in \cite{Poirier} was interested to
construct the graded dual to the Mantaci-Reutenauer Hopf algebra. His goal was
to produce analogues of $QSym$ and $NSym$ Hopf algebras for type $B$. For
this, Poirier introduced a larger algebra containing the desired Hopf
algebras. This larger algebra is isomorphic to our $DQSym$. More
generalizations and properties of these Hopf algebras are studied in
\cite{Hohlweg, NovThi}.
Moreover, in an upcoming work, S. Kapetanovich \cite{Kape} shows that
$\DQSym$ is the terminal object in the category of bigraded Hopf algebras
equipped with a multiplicative character. She also computes explicitly its
even and odd Hopf  subalgebras
\end{proof}

\section{Conjectures}

\subsection{The quotient $\Q[{\bf x},{\bf y}]/\langle \DQSym_n^+\rangle$.}We are now interested in the study of the ideal generated by diagonally quasi-symmetric polynomials without constant term, \ie 
$$\DI_n=\langle M_{\mathfrak a}({\bf x},{\bf y})\rangle_{{\mathfrak a}\neq\emptyset}.$$
This ideal is clearly bihomogeneous.
Next, in the spirit of our discussion in section \ref{intro}, we consider the quotient space
    $$\DQ_n:=\Q[{\bf x},{\bf y}]/\DI_n,$$ 
as well as the space of {\em diagonal $TL_n$-harmonic polynomials} as 
   $$\DSH_n:=\DI_n^\perp$$ 
 with respect to the scalar product of \pref{scalXY}. These two spaces are isomorphic as bigraded modules. 
Observe that
\begin{equation}\label{tlhxy}
     \TLH_n\simeq \DSH_n\cap \Q[{\bf x}], \qquad {\rm and}\qquad \TLH_n\simeq \DSH_n\cap\Q[{\bf y}].
 \end{equation}
Recall that the {\em bigraded Hilbert series} of a bigraded vector subspace $\V$   is by definition
    $$H_\V(q,t)=\sum_{h,k}q^ht^k\dim(\V_{h,k})$$
where $\V_{h,k}$'s are the components appearing in the decomposition
   $$\V=\bigoplus_{h,k} \V_{h,k},$$
 corresponding to the bigraded structure of $\V$. The Hilbert series considered here are in fact polynomials.
To stress out the symmetries of these polynomials, we use a cartesian coordinates entries matrix notation for them. 
Thus, for a polynomial $p(q,t)$  we define the corresponding {\em Hilbert matrix}, $M$, to have entries
    $$M(i,j)={\rm coefficient\ of\ } t^{k-i+1 }q^{j-1}\ {\rm in\ } p(q,t),$$
with $k$ equal to the degree in $t$ of $p(q,t)$.
Thus we have the correspondence
$$\begin{array}{l}
  3\,t^2+\\
  2\,t+3\,tq+\\
   1+3\,q+4\,q^2+q^3
 \end{array}\quad \leftrightarrow \qquad
\begin{pmatrix}
3&&&\\
2&3&&\\
1&3&4&1\\
\end{pmatrix}$$
Here, we omit $0$ entries.

\begin{conj}\label{conjmat}
The Hilbert matrix, of the space $\DQ_n$, is the $n\times n$ lower triangular matrix, $M_n$, with all entries in the main diagonal equal to the $(n-1)^{\rm th}$ Catalan number:
   $$M_n(i,i)=\frac1n{2(n-1)\choose(n-1)},$$
and all other non zero entries recursively obtained as:
\begin{equation}\label{conj_rec}
   M_n(i,j)=\sum_{\twolines{i'\geq i}{j'\leq j}}M_{n-1}(i',j'),
\end{equation}
when $i>j$.
\end{conj}
\noindent For instance, the Hilbert matrices $M_n$, for $1\le n\le5$, are 
\begin{eqnarray*}
&&M_1=\begin{pmatrix}
1
\end{pmatrix},
\qquad
M_2=\begin{pmatrix}
1&\\
1&1
\end{pmatrix},
\qquad
M_3=\begin{pmatrix}
2&&\\
2&2&\\
1&2&2
\end{pmatrix},
\\
&&M_4=\begin{pmatrix}
5&&&\\
5&5&&\\
3&7&5&\\
1&3&5&5
\end{pmatrix},
\quad
M_5=\begin{pmatrix}
14&&&&\\
14&14&&&\\
9&24&14&&\\
4&14&24&14&\\
1&4&9&14&14
\end{pmatrix}.
\end{eqnarray*}
Actually, the values in the first column (and last line) of $M_n$ are known to satisfy recurrence (\ref{conj_rec}). In view of (\ref{tlhxy}) and the isomorphism of graded spaces $\QQ_n\simeq \TLH_n$, this follows immediately from our result of \cite{ABBQ}, giving the following explicit formula for the Hilbert series of $\QQ_n$:
\begin{equation}
    H_{\QQ_n}(q)=\sum_{k=0}^{n-1}\frac{n-k}{n+k}{n+k \choose k} q^k\,.
\end{equation}
In \cite{ABBQ}, we have in fact constructed a nice explicit linear basis of monomials for $\QQ_n$. We are going to describe here, for each $n$, a set $\B_n$ of monomials which we conjecture to be a linear basis of $\DQ_n$. As should be expected this set contains two copies of the basis for $\QQ_n$; one in the ${\bf x}$-variables, and one in the ${\bf y}$-variables. The degrees of monomials in $\B_n$ are all less or equal to $n-1$. Moreover, the construction will make clear that  the bigraded enumeration of $\B_n$ agrees with the Hilbert matrix of Conjecture \ref{conjmat}. First, whenever $i+j<n-1$, we construct recursively the subset $\B_n(i,j)$ of $\B_n$, consisting of its elements of bidegree $(i,j)$, by setting
\begin{equation}\label{base_rec}
   \B_n(i,j):=\bigcup_{\twolines{a\le i}{b\le j}}\B_{n-1}(a,b)\cdot x_n^{i-a}y_n^{j-b}.
 \end{equation}
Here, in the right hand side, the product  is to be understood as
   $$B\cdot x_n^cy_n^d:=\{{\bf z}\cdot  x_n^cy_n^d\ |\ {\bf z}\in B \}.$$
Second, starting with $\B_n(n-1,0)$, also constructed via recurrence (\ref{base_rec}), we generate successively the sets  $\B_n(n-i-1,i)$, for $i$ running from $2$ to $n-1$:
 \begin{equation}\label{base_diag}
     \B_n(n-i-1,i):=\{ \Phi({\bf z})\ |\ {\bf z}\in \B_n(n-i,i-1)\ \}
 \end{equation}
with $\Phi$ being defined as in (\ref{phi}).
Consider the operation $\phi$, defined on bicompositions  as $\phi({\bf a},{\bf b}):=({\bf a}-{\bf e}_m,{\bf b}+{\bf e}_m)$,
 with ${\bf e}_m$ the usual $m^{\rm th}$ unit vector, and $m$ equal to the smallest index for which ${\bf a}$ has a non zero coordinate. We then define the operation $\Phi$ on monomials, by simply setting
\begin{equation}\label{phi}
    \Phi(({\bf x}{\bf y})^\mathfrak{a}_{\bf i}):=({\bf x}{\bf y})^{\phi(\mathfrak{a})}_{\bf i}.
 \end{equation}

\begin{conj}\label{conjbas}
  The set $\B_n$, constructed using rules {\rm (\ref{base_rec})} and {\rm (\ref{base_diag})},  is a linear basis of $\DQ_n$.
\end{conj}

For $1\le n\le3$, the bases generated by the recursive rule above are
\begin{eqnarray*}
\B_1&=&\{1\}\\
\B_2&=&\left\{\begin{matrix} y_2&\\1&x_2\end{matrix}\right\}\\
\B_3&=& \left\{\begin{matrix}y_3^2,y_2y_3\\ y_3,y_2&y_3x_3,y_2x_3\\1&x_3,x_2 &x_3^2,x_2x_3\end{matrix}\right\}
\end{eqnarray*}
They are presented here in matrix form to stress out how this construction is in agreement with the Conjecture \ref{conjmat}.

\subsection{The quotient $\DQSym_n/\langle \DSym_n^+\rangle$.}\label{qsym_sym}
Since we can check that 
 $$\DSym=\Q[M_{\alpha}\ |\ \alpha\in \N^2],$$
the quotient of $\DQSym$, by the ideal generated by constant term free elements of $\DSym$, is easy to characterize. Moreover, we clearly have a bigraded module isomorphism 
   $$\DQSym\simeq \DSym\otimes \DQSym/\langle \DSym^+\rangle.$$
It follows that $\DQSym$ is a free $\DSym$ module. Recall that this is in the context when the number of variables has been sent to infinity. The story is not so simple when we restrict to $n$ the number of ${\bf x}$ and ${\bf y}$ variables. As we will see, we lose the analogous freeness property in this bivariate polynomials case. This is in contrast with the univariate case for which we have indeed
\begin{equation}\label{univarie}
  \QSym_n\simeq \Sym_n\otimes \QSym_n/\langle \Sym_n^+\rangle,
\end{equation}
as conjectured in \cite{BerReu}, and then shown in \cite{GW}. It readily follows from  (\ref{univarie}) that we have the Hilbert series formula
 \begin{equation}\label{formule_R}
    H_{{\rm R}_n}(q)=(1+q)\,\cdots\,(1+q+\ldots+q^{n-1})\,\frac{(1-q)^{n+1}-q^{n+1}}{1-2\,q},
  \end{equation}
 where we have denoted ${\rm R}_n$ the space  $\QSym_n/\langle \Sym_n^+\rangle$.
 It is not immediate that the polynomials in right hand side of (\ref{formule_R}) have indeed positive integer coefficients. This was shown to be true in \cite{BerReu}, using the fact that these polynomials are solutions of the recurrence
        $$\Psi_n(q)=\Psi_{n-1}(q)+q^n(n!_q-\Psi_{n-1}(q)),\qquad n\geq 1$$
with initial condition $\Psi_0(q)=1$. The first values of these $\Psi_n(q)$ are
   \begin{eqnarray*}
        \Psi_1(q)&=&1\\
       \Psi_2(q)&=&1+q^3\\
       \Psi_3(q) &=& 1+q^3+2\,q^4+2\,q^5\\
   \end{eqnarray*}
We want to underline that $\DQSym_n$ cannot be a free $\DSym_n$ module. Otherwise, the Hilbert series of the space
    $$\mathcal{R}_n:= \DQSym_n/\langle \DSym_n^+\rangle,$$
 would have to be equal to
 \begin{equation}\label{mauvaise}
    \frac{\displaystyle \sum_{k=0}^n ((1-t)^{-1}(1-q)^{-1}-1)^k}{h_n[(1-t)^{-1}(1-q)^{-1}]},
 \end{equation}
 using plethystic substitution notation. Recall that this means that
    $$h_n\left[\frac{1}{(1-t)(1-q)}\right]=\sum_{\lambda\vdash n}
           \frac{1}{z_\lambda} \prod_{k=1}^{\ell(\lambda)} \frac{1}{(1-q^k)(1-t^k)},$$
 with $z_\lambda=\prod_k k^{d_k} d_k!$, if the partition $\lambda$ has $d_k$ parts of size $k$.
 One can readily check that (\ref{mauvaise}) contains some negative integer coefficients in general. This is so even in the simple case of $n=2$ for which we get:
   $$\frac{(1+t+q-q\,t)(1-t^2)(1-q^2)}{1+q\,t}=1+qt+{q}^{3}+{q}^{2}t+q{t}^{2}+{t}^{3}-{q}^{3}t-{q}^{2}{t}^{2}-q{t}^3+\ldots$$
Hence, the Hilbert series cannot be given by (\ref{mauvaise}).
The actual values for the Hilbert series of $\mathcal{R}_n$ are given below (in matrix form), for $n\leq 3$.
$$
\begin{pmatrix} 1\end{pmatrix},\qquad
\begin{pmatrix} 
1\\
0&1\\
0&1&1\\
1&0&0&1
\end{pmatrix},\qquad
 \begin{pmatrix} 
2\\
2&2\\
1&4&2\\
0&3&4&2\\
0&1&3&4&2\\
1&0&0&1&2&2
\end{pmatrix}.
$$
Clearly, the first column (and first line) correspond to instances of (\ref{formule_R}). Further computations suggest many other nice properties of these Hilbert series, with proofs using the techniques of \cite{BerReu} extended to this new context. This will be the subject of an upcoming paper.



\begin{thebibliography}{10}   

\bibitem{artin}   {\sc E. Artin},      {\em Galois Theory},  Notre Dame Mathematical Lecture {\bf 2} (1944), Notre Dame, IN.  

\bibitem{ABS}   {\sc M. Aguiar, N. Bergeron and F. Sottile},   {\em  Combinatorial Hopf Algebras and Generalized Dehn-Sommerville Relations}, to appear. ArXive math.CO/0310016.

\bibitem{AS}   {\sc M. Aguiar and F. Sottile},    {\em Structure of the Malvenuto-Reutenauer Hopf algebra of permutations}, to appear. ArXive math.CO/0203282.

\bibitem{a9}     {\sc J.-C. Aval and N. Bergeron}, {\em Catalan Paths and Quasi-Symmetric Functions}, Proc. Amer. Math. Soc.,  {\bf 131} (2003) 1053--1062. 

\bibitem{ABBQ}   {\sc J.-C. Aval, F. Bergeron and N. Bergeron},   {\em Ideals of Quasi-Symmetric Functions and Super-Covariant Polynomials for $\S_n$},  Adv. Math, {\bf 181-2} (2004)  353-367.

 \bibitem{lattice}   {\sc F. Bergeron, N. Bergeron, A. Garsia, M. Haiman and G. Tesler}, {\em  Lattice Diagram Polynomials and Extended Pieri Rules}, Adv. Math., {\bf 142}   (1999), 244-334.   
 
\bibitem{nabla}
  {\sc F. Bergeron, A. Garsia, M. Haiman and G. Tesler}, {\em
  Identities and Positivity Conjectures for Some Remarkable Operators in the Theory of Symmetric Functions}, Methods and Applications of Analysis, Volume 6, No. 3, (1999), 363--420.

\bibitem{Berg2}   {\sc F. Bergeron, A. Garsia and G. Tesler}, {\em Multiple Left Regular  Representations Generated by Alternants}, J. of Comb. Th., Series A, {\bf 91,   1-2}  (2000), 49--83.  

\bibitem{BerReu}   {\sc F. Bergeron, and C. Reutenauer}, {\em The coinvariant space of quasisymmetric polynomials}, in preparation.  

\bibitem{BMSW}   {\sc N.~Bergeron, S.~Mykytiuk, F.~Sottile, and S.~van Willigenburg}, {\em     Pieri Operations on Posets},  \newblock J. of Comb. Theory, Series A, {\bf 91} (2000), 84--110 .  

\bibitem{borel} {\sc A. Borel}, {\em Sur la cohomologie des espaces fibr\'es principaux et des espaces homog\`enes des groupes de Lie compacts}, Ann. of Math., {\bf 57}  (1953), 115--207.  

\bibitem{CP}   {\sc C. de Concini and C. Procesi}, {\em Symmetric functions, conjugacy  classes and the flag variety},  Invent. Math. {\bf 64} (1981), 203--230.  

\bibitem{CLO}   {\sc D. Cox, J. Little and D. O'Shea}, {\em Ideals, Varieties, and Algorithms},   Springer-Verlag, New-York, 1992.  

\bibitem{GH}   {\sc A.~M. Garsia and M.~Haiman}, {\em A graded representation model for     {M}acdonald's polynomials}, Proc. Nat. Acad. Sci. U.S.A. {\bf 90} (1993),      no.~8, 3607--3610.  

\bibitem{orbit} A. Garsia and M. Haiman, {\it Orbit Harmonics and Graded Representations}, \'Editions du Lacim, to appear.  

\bibitem{GW} A. Garsia and N. Wallach {\it }, {\em  Qsym over Sym is free.} J. Combin. Theory Ser. A {\bf 104} (2003), no. 2, 217--263.

\bibitem{NC}     {\sc I.~Gelfand, D.~Krob, A.~Lascoux, B.~Leclerc, V.~Retakh, and J.-Y.  Thibon},      {\em Noncommutative symmetric functions}, Adv. in Math., {\bf 112} (1995),     ~218--348.   

\bibitem{ges}   {\sc I.~Gessel}, {\em Multipartite ${P}$-partitions and products of skew     {S}chur functions}, in Combinatorics and Algebra (Boulder, Colo., 1983),     C.~Greene, ed., vol.~34 of Contemp. Math., AMS, 1984, pp.~289--317.

 \bibitem{haiman}     {\sc M. Haiman}, {\em {H}ilbert schemes, polygraphs, and the Macdonald  positivity conjecture}, J. Amer. Math. Soc., {\bf 14} (2001), 941-1006.  
 
\bibitem{hivert}     {\sc F. Hivert}, {\em Hecke algebras, difference operators, and  quasi-symmetric functions\/}, Adv. in Math., {\bf 155} (2000),  181--238.   

\bibitem{Hoff}
  {\sc M. E. Hoffman}, {\em  Quasi-shuffle products}, J. Algebraic Combin. {\bf 11} (2000), no. 1, 49--68.
   \bibitem{jones} {\sc V. Jones}, {\em Index for subfactors}, Inv. Math. {\bf 72} (1983) 1--25.  
   
\bibitem{Hohlweg}
{\sc    P. Baumann and C. Hohlweg}, {\em A Solomon theory for the wreath
products $G\wr\mathfrak S_n$}.  Preprint (2004)

\bibitem{Kape}
   {\sc    S. Kapetanovich}, {\em Dissertation Thesis York}. In preparation. 
   

\bibitem{NCSF4}{\sc D. Krob} and {\sc J.-Y. Thibon}, {\it Noncommutative
    symmetric functions IV: Quantum linear groups and hecke algebras at
    $q=0$}, J. Alg. Comb. 
    {\bf 6} (1997), 339i--376.
%
 
\bibitem{mac}   {\sc I.~Macdonald}, {\em Symmetric Functions and Hall Polynomials}, Oxford     Univ.~Press, 1995,  \newblock second edition.  

\bibitem{MR}  {\sc C. Malvenuto and C. Reutenauer}, {\em Duality     between quasi-symmetric functions     and the {S}olomon descent algebra},     Journal of  Algebra, \textbf{177} (1995), 967--982. 
   
\bibitem{MM65} {\sc J.~W. Milnor and J.~C. Moore}, \emph{On the structure of {H}opf algebras},    Ann. of Math.  81 (1965), no.~2, 211--264.    

\bibitem{NovThi}
{\sc J.-C. Novelli} and {\sc J.-Y.  Thibon}, {\em Free Quasi-symmetric Functions of arbitrary level}.  To appear. math.co/0405597.

\bibitem{Poirier}
{\sc    S. Poirier}, {\em Cycle type and descent set in wreath products}.
Discrete mathematics {\bf 180} (1998) p. 315--343.

\bibitem{stanley1}   {\sc R.~Stanley}, {\em Enumerative Combinatorics, Vol.~1}, Wadsworth and     Brooks/Cole, 1986.  

\bibitem{stanley2}     {\sc R.~Stanley}, {\em Enumerative Combinatorics Vol.~2}, no.~62 in  Cambridge      Studies in Advanced Mathematics, Cambridge University Press, 1999.  \newblock Appendix 1 by Sergey Fomin.  

\bibitem{steinberg}    {\sc R.~Steinberg}, {\em Differential equations invariant under  finite  reflection groups}, Trans. Amer. Math. Soc., {\bf 112} (1964), 392--400.

%
\bibitem{Stem}
 {\sc J. Stembridge}, {\em Enriched $P$-partitions}, Trans. Amer. Math. Soc.
{\bf 349} (1997), 763--788.

  
\end{thebibliography}
\end{document}